# Arbitrary order finite volume well-balanced schemes for the Euler equations with gravity

C. Klingenberg*     G. Puppo†     M. Semplice‡

July 6, 2018

**Abstract**

This work presents arbitrary high order well balanced finite volume schemes for the Euler equations with a prescribed gravitational field. It is assumed that the desired equilibrium solution is known, and we construct a scheme which is exactly well balanced for that particular equilibrium. The scheme is based on high order reconstructions of the fluctuations from equilibrium of density, momentum and pressure, and on a well balanced integration of the source terms, while no assumptions are needed on the numerical flux, beside consistency. This technique allows to construct well balanced methods also for a class of moving equilibria. Several numerical tests demonstrate the performance of the scheme on different scenarios, from equilibrium solutions to non steady problems involving shocks. The numerical tests are carried out with methods up to fifth order in one dimension, and third order accuracy in 2D.

## 1 Introduction

In this paper we are concerned with the numerical approximation of the flow of a gas in a gravitational field. The problem is modelled with Euler gas dynamics equations with a source term, containing the gravitational force, namely

$$\begin{cases} \dfrac{\partial \rho}{\partial t} + \nabla \cdot (\rho v) &= 0 \\ \dfrac{\partial (\rho v)}{\partial t} + \nabla \cdot (\rho v \otimes v + pI) &= -\rho \nabla \Phi \\ \dfrac{\partial E}{\partial t} + \nabla \cdot (v(E+p)) &= -\rho v \nabla \Phi \end{cases} \quad (1)$$

*University of Würzburg, Würzburg, Germany (*klingen@mathematik.uni-wuerzburg.de*)
†Università dell'Insubria, Como, Italy (*gabriella.puppo@uninsubria.it*)
‡Università di Torino, Torino, Italy (*matteo.semplice@unito.it*)

for $x \in \mathbb{R}^d$, with $\rho \geq 0$ being the density, $v \in \mathbb{R}^d$ the velocity, $m = \rho v$ the momentum, $p \geq 0$ the pressure and $E$ the total energy per unit volume. Further, the internal energy density is $e$ and it is given by $\rho e = E - \frac{1}{2}\rho v \cdot v$. Pressure is determined from $e$ and $\rho$ through the equation of state (EOS). For an ideal gas, the internal energy depends only on the temperature $e = e(T)$, but other cases are possible. The state of a gas is determined by only two thermodynamic variables, as, for instance, the pressure and the density. Thus, pressure, temperature and density determine a triplet of functions, in which any one of them is determined by the remaining two through the equation of state. For example, for an ideal gas, $p = \rho RT$, where $R$ is the universal constant for an ideal gas.

We will suppose that system (1) is completed with an initial condition

$$\rho(x, t=0) = \rho_0(x), \qquad v(x, t=0) = v_0(x), \qquad p(x, t=0) = p_0(x). \qquad (2)$$

The numerical integration of equations of the form (1) presents several challenges: singularities may form in a finite time, even from smooth initial data, making the solution rich in structure. For this reason, high order accurate schemes, tailored to deal with discontinuities, are particularly interesting. They permit to resolve fine scales on the solution even using relatively coarse grids. For a classical review of the issues relevant in the construction of high order non oscillatory schemes, see [21].

The main focus of this paper is on a further challenge in the integration of balance laws, which is due to the presence of the source term. System (1) can be endowed with non trivial steady states, and often small perturbations of such equilibrium states are of particular interest in applications. However, small perturbations of steady states may go totally undetected, if they are of the same size of the local truncation error. For this reason, much research has concentrated on the development of *well balanced schemes*, which are able to preserve steady states exactly at the discrete level, thus enabling also the detection of small perturbations of steady states.

The development of well balanced schemes started on the system of the shallow water equations, with, initially, the goal to preserve the lake at rest solution. Pioneering works in this field are [2, 14], but the literature on this topic is huge. Here we mention especially the high order well balanced schemes of [18, 19] and the technique of well balancing thanks to a hydrostatic reconstruction of [1], which employ tools that are at the basis of our approach to well balancing under a gravitational field. See also the extension proposed in [24] for the high order preservation of moving water equilibria.

More recently, new applications of well balanced schemes have been considered. In particular, gas dynamics flows in a gravitational field are endowed with non trivial equilibria which are of particular interest for astrophysical and metheorological applications. In this work, we concentrate on the construction of high order numerical methods which are well balanced for equilibrium solutions of Euler equations with gravity. A pioneering work in this field is [4], which constructs a well-balanced scheme transforming the source terms into numerical fluxes at equilibrium. Many recent papers deal with this problem. Initially, well balanced schemes were written for particular equilibria, as in [25], where high order finite difference methods for isothermal equilibrium are proposed, and [13], which concentrates on isentropic flows, for general equations of state. The work [5] proposes a second order scheme which is well balanced thanks to an auxiliary function, which determines an ad hoc recontruction. In [11, 23], well balancing is achieved with a relaxation scheme, which includes the enforcement of a steady state within the numerical flux. We also mention [6, 17] which are concerned with Discontinuos Galerkin methods. The work [3] concentrates on well balancing on low Mach smooth steady states.

A few very recent schemes, [7] or [10], are able to detect equilibrium states automatically, and be well balanced only against such states. These schemes are only second order accurate, and they are well balanced with respect to an approximation of the exact unknown steady state.

The method we propose in this work assumes that a particular equilibrium is given, around which the scheme is well balanced, as in [16]. Our method will preserve this steady state exactly, and will be able to resolve accurately very small perturbations around it. In many applications, it is reasonable to assume that one knows the structure of the steady state of interest: for instance, one may expect that at equilibrium the flow will be isentropic, or isothermal. Our method is built on the idea that the hydrostatic equilibrium around which the scheme is well balanced is chosen by the user. This allows us to construct a scheme which is well balanced also for a type of moving equilibrium, which, as far as we know, has not been considered yet in the literature on well balanced schemes for Euler with gravity. We believe that this new equilibrium can be of interest for applications, because it permits to compute high order well balanced solutions also for the case of the uniform flow of gas perpendicularly to a gravitational field.

The scheme we propose is built on three main ideas: we propose to reconstruct variables as fluctuations from the prescribed steady state, because any piecewise polynomial reconstruction algorithm is able to reproduce constants exactly. This eliminates artefacts due to the artificial diffusion inherent in the stabilization terms, see also [1]. Secondly, the prescribed equilibrium is used to eliminate the gravitational potential from the equations, which is used to construct a well balanced second order accurate quadrature for the source. Finally, we use extrapolation, as in [18] to boost the second order well balanced quadrature to any desired order.

The rest of the paper is organized as follows. We start with a discussion on the numerical treatment of Euler with gravity steady state solutions in §2. The numerical scheme and its properties are discussed in Sections 3 for the one dimensional case and 4 for higher space dimensions. The results of the numerical tests are then reported in Section 5 and our conclusions are summarized in Section 6.

## 2   Steady state solutions and well balanced schemes

System (1) is a balance law, and thus it is possible to have non trivial steady states, when the flux terms balance exactly the source. In particular, Euler equations with gravity can have hydrostatic equilibrium solutions, if the following conditions are satisfied

$$v(x,t) \equiv 0 \qquad \nabla p = -\rho \nabla \Phi. \tag{3}$$

Since system (1) does not have dissipation terms, steady states can occur only if they are present already at the level of the initial condition, and are consistent with boundary conditions, if present. Not all initial conditions are consistent with a possible steady state, since the existence of a hydrostatic steady state requires that

$$\nabla \times (\rho \nabla \Phi) = 0. \tag{4}$$

In this work we will suppose that, given the gravitational potential $\Phi$, two scalar functions $\alpha$ and $\beta$ are known such that

$$\nabla \beta = -\alpha \nabla \Phi. \tag{5}$$

So, if at some time $\rho(x,t) = \alpha(x)$, the consistency condition (4) is satisfied and an hydrostatic steady state becomes possible. In many cases, one is interested in preserving a particular equilibrium state, from which the functions $\alpha$ and $\beta$ can be derived *a-priori*. Typical cases include, but are not limited to, the following ones.

If the temperature is constant, equilibrium for an ideal gas is described by

$$\alpha^{\text{iso}}(x) = \frac{e^{-\Phi(x)/T_{eq}}}{T_{eq}} \qquad \beta^{\text{iso}}(x) = e^{-\Phi(x)/T_{eq}}. \tag{6}$$

In fact, for an ideal gas, $p = \rho T$ (we are taking the gas constant $R = 1$ for simplicity). Thus the equilibrium equation becomes $\nabla p = -p/T \, \nabla \Phi$. Since the temperature is constant $T \equiv T_{eq}$, the equation becomes $\nabla(\ln p) = -\nabla(\Phi/T)$, and one can integrate both sides, obtaining the same result irrespective of the path.

Similarly, if the desired equilibrium targets a polytropic gas, $p\rho^{-\nu} = C$ (constant), the functions $\alpha$ and $\beta$ are given by

$$\alpha^{\text{poly}}(x) = \left(1 - \tfrac{\nu-1}{\nu}\Phi(x)\right)^{\frac{1}{\nu-1}} \qquad \beta^{\text{poly}}(x) = (\rho(x))^{\nu}. \tag{7}$$

Note that this in particular includes the isentropic equilibrium when $\nu = \gamma$.

These are classical solutions of clear interest for physical applications. But this does not by any means cover all possible cases. For instance, if the density is constant, we would take $\beta = -\rho\Phi$. Further, in the numerical tests, we will also consider the following equilibrium around a potential with radial symmetry,

$$\alpha^{\text{gen}}(r,\theta) = e^{-r} \qquad \beta^{\text{gen}}(r,\theta) = (1+r)e^{-r}, \qquad \Phi(r) = r^2. \tag{8}$$

Once the functions $\alpha$ and $\beta$ are known, if $v_0(x,y) = 0$, $p_0 = K\beta$ and $\rho_0 = K\alpha$, for some constant $K$, then the solution will remain stationary for all time.

Another interesting equilibrium which so far, to the best of our knowledge, has not been studied, concerns equilibrium solutions with non zero constant speed in multi dimensions, as in the case of a steady breeze along an horizontal surface. We consider a gravitational field with a constant direction, and we align the system of reference with the $y$ axis parallel to $\nabla \Phi$. Then a solution of the form

$$\rho(x,y) = \alpha(y), \quad p(x,y) = \beta(y), \quad v(x,y) = (U, 0), \tag{9}$$

with $\alpha$ and $\beta$ satisfying (5) is a steady state solution, for any constant $U$. The well balanced discretization of this equilibrium is introduced in §4.1.

A standard discretization of system (1) in general fails to preserve steady states exactly. There are two main issues at stake. We illustrate the origin of the failure to preserve steady states considering a simple first order finite volume scheme for (1). For simplicity, we will consider a uniform grid in space. The computational domain is covered with control volumes $V_j^n = (x_j - \frac{\Delta x}{2}, x_j + \frac{\Delta x}{2}) \times (t^n, t^n + \Delta t)$, where $\Delta x$ and $\Delta t$ are the grid spacings in space and time respectively, $\lambda = \Delta t/\Delta x$ and $x_j = j\Delta x, j \in \mathbb{Z}$, $t^n = n\Delta t, n \in \mathbb{N}$. A first order discretization of system (1) will give

$$\overline{U}_j^{n+1} = \overline{U}_j^n - \lambda\left[\mathcal{F}\left(\overline{U}_{j+1}^n, \overline{U}_j^n\right) - \mathcal{F}\left(\overline{U}_j^n, \overline{U}_{j-1}^n\right)\right] + \lambda S(\overline{U}_j^n).$$

Here $\mathcal{F}(a,b)$ is a standard numerical flux, given by

$$\mathcal{F}(a,b) = \tfrac{1}{2}(f(a) + f(b)) - \tfrac{1}{2}Q(a,b)(a-b),$$

where $Q(a,b)$ is the viscosity matrix of the numerical method. For instance, $Q(a,b) = \mu \mathbb{I}$ for the Lax Friedrichs numerical flux, where $\mu$ is the artificial diffusion coefficient, or $Q(a,b)$ is the Roe matrix for the Roe numerical flux. If we specialize this discretization to a steady state solution of the form (3), for the Lax Friedrichs numerical flux, we find

$$\rho_j^{n+1} = \rho_j^n + \tfrac{1}{2}\mu\lambda(\rho_{j+1}^n - 2\rho_j^n + \rho_{j-1}^n)$$
$$(\rho v)_j^{n+1} = -\tfrac{1}{2}\lambda(p_{j+1}^n - p_{j-1}^n) - \rho_j^n\left(\Phi(x_{j+1/2}) - \Phi(x_{j-1/2})\right)$$
$$E_j^{n+1} = E_j^n + \tfrac{1}{2}\mu\lambda(E_{j+1}^n - 2E_j^n + E_{j-1}^n).$$

The second equation generates momentum spuriously because the discretization of the source does not match the differences in the pressure. The remaining equations move the density and the energy away from the steady state, because of the artificial diffusion term. Thus, to achieve well balancing, two aspects should be considered

- reconstruct along equilibrium variables, to ensure that at equilibrium the two interface states on which the numerical flux is built coincide: then the consistency of the numerical flux implies that $\mathcal{F}(U,U) = f(U)$, with no artificial diffusion;

- write a well balanced quadrature of the source, to ensure that the numerical flux and the cell average of the source balance exactly at the discrete level, for equilibrium solutions.

Following this framework, it is possible to write well balanced numerical schemes for any order of accuracy, and for any consistent numerical flux.

## 3 One-dimensional numerical scheme

Let us first describe the discretization of system (1) in one space dimension. We suppose we are given two functions $\alpha(x), \beta(x)$, as in (5), that is with $\nabla \beta = -\alpha \nabla \Phi$. Then, whenever $u = 0$, $p(x,t) = \beta(x)$ and $\rho(x,t) = \alpha(x)$, the system is in hydrostatic equilibrium, and the solution will remain constant.

The idea is to reconstruct the *fluctuations* of $\rho(x,t)$ and $p(x,t)$ from $\alpha$ and $\beta$, so that at equilibrium the reconstruction is identically zero.

Since $\alpha$ and $\beta$ are known functions, we can compute their cell averages, $\overline{\alpha}_j$ and $\overline{\beta}_j$. For a fixed time $t$, introduce auxiliary variables

$$r(x,t) = \rho(x,t) - \alpha(x) \qquad \pi(x,t) = p(x,t) - \beta(x). \tag{10}$$

Since, during the reconstruction, the time is fixed, we temporarily drop the $t$ dependence. We start describing the first and second order schemes. Up to second order we can take

$$\overline{p} = (\gamma - 1)(\overline{E} - \tfrac{1}{2}\overline{m}^2/\overline{\rho}). \tag{11}$$

Then, we can compute the cell averages of the auxiliary variables,

$$\overline{r}_j = \overline{\rho}_j - \overline{\alpha}_j \qquad \overline{\pi}_j = \overline{p}_j - \overline{\beta}_j.$$

We apply a non oscillatory reconstruction to $\overline{r}$, $\overline{\pi}$ and $\overline{m} = \overline{\rho v}$. Let $r^{\pm}_{j+1/2}$ and $\pi^{\pm}_{j+1/2}$ be the left and right reconstructed values at the cell interface located in $x = x_j + 1/2 h$ for the density and pressure fluctuations. Then the reconstructed values for the density and the pressure can be recovered as

$$\rho^{\pm}_{j+1/2} = r^{\pm}_{j+1/2} + \alpha(x_{j+1/2}) \tag{12a}$$

and

$$p^{\pm}_{j+1/2} = \pi^{\pm}_{j+1/2} + \beta(x_{j+1/2}). \tag{12b}$$

Note that the accuracy of the reconstructed data $\rho^{\pm}_{j+1/2}$ and $p^{\pm}_{j+1/2}$ is of the same order $q$ one would achieve reconstructing the point values directly from $\overline{\rho}$ and $\overline{p}$, with the same reconstruction. This is due to the continuity of the functions $\alpha$ and $\beta$ across the interfaces. The reconstructed values for the momentum do not need to be modified, because $(\rho v)^{\pm}_{j+1/2} = 0$ at equilibrium. Then we recover the energy as

$$E^{\pm}_{j+1/2} = \tfrac{1}{2}(\rho v)^{\pm}_{j+1/2}{}^2/\rho^{\pm}_{j+1/2} + p^{\pm}_{j+1/2}/(\gamma - 1) \tag{12c}$$

As we will see below, at equilibrium, the reconstructed data are continuous, i.e. $U^+_{j+1/2} = U^-_{j+1/2} = U_{j+1/2}$. Then the numerical flux is

$$\mathcal{F}\left(U^+_{j+1/2}, U^-_{j+1/2}\right) = \mathcal{F}\left(U_{j+1/2}, U_{j+1/2}\right) = f(U_{j+1/2}),$$

and no artificial diffusion is introduced. The last task is to discretize the source. To this end, note that the source can be written exactly as

$$-\rho \nabla \Phi = -\frac{\rho}{\alpha} \alpha \nabla \Phi = \frac{\rho}{\alpha} \nabla \beta. \tag{13}$$

In our scheme, we use this formula to discretize the source term. In particular, the quadrature $Q_j^{\rho v}$, giving the cell average of the source in the momentum equation, is defined as

$$Q_j^{\rho v} = \frac{1}{2} \left( \frac{\rho^+_{j-1/2}}{\alpha(x_{j-1/2})} + \frac{\rho^-_{j+1/2}}{\alpha(x_{j+1/2})} \right) \frac{\beta(x_{j+1/2}) - \beta(x_{j-1/2})}{\Delta x}. \tag{14}$$

Similarly, the quadrature rule for the source in the energy equation is

$$Q_j^E = \frac{1}{2} \left( \frac{(\rho v)^+_{j-1/2}}{\alpha(x_{j-1/2})} + \frac{(\rho v)^-_{j+1/2}}{\alpha(x_{j+1/2})} \right) \frac{\beta(x_{j+1/2}) - \beta(x_{j-1/2})}{\Delta x}. \tag{15}$$

These two formulas yield a second order accurate approximation of $\overline{S(U)}_j$.

In general, we will define the well balanced quadrature for a function $y$ with the expression

$$Q_j^y = \frac{1}{2} \left( \frac{y^+_{j-1/2}}{\alpha(x_{j-1/2})} + \frac{y^-_{j+1/2}}{\alpha(x_{j+1/2})} \right) \frac{\beta(x_{j+1/2}) - \beta(x_{j-1/2})}{\Delta x}. \tag{16}$$

**Theorem 1.** *The semidiscrete scheme*

$$\frac{d}{dt}\overline{U}_j(t) = -\frac{1}{\Delta x}\left[\mathcal{F}\left(U^+_{j+1/2}(t), U^-_{j+1/2}(t)\right) - \mathcal{F}\left(U^+_{j-1/2}(t), U^-_{j-1/2}(t)\right)\right] + Q_j(t), \tag{17}$$

*where $Q_j = [0, Q_j^{\rho v}, Q_j^E]$ is computed in (14) and (15), while the boundary extrapolated data are obtained with the well-balanced reconstructions (12), is exactly well balanced on the hydrostatic equilibrium solution $\rho(x,t) = \alpha(x)$ and $p(x,t) = \beta(x)$.*

*Proof.* If the data at time $t$ are in hydrostatic equilibrium, $\overline{\rho v}$, $\overline{r}_j$ and $\overline{\pi}_j$ are zero for all $j$. Then, the boundary extrapolated data for the pressure are $p^\pm_{j+1/2} = \beta(x_{j+1/2})$. The momentum flux reduces to $\mathcal{F}^{\rho v}_{j+1/2} = \beta(x_{j+1/2})$, and the equation for momentum becomes

$$\frac{d}{dt}\overline{\rho v}_j(t) = -\frac{\beta(x_{j+1/2}) - \beta(x_{j-1/2})}{\Delta x}$$
$$+ \frac{1}{2}\left(\frac{\rho^+_{j-1/2}}{\alpha(x_{j-1/2})} + \frac{\rho^-_{j+1/2}}{\alpha(x_{j+1/2})}\right)\frac{\beta(x_{j+1/2}) - \beta(x_{j-1/2})}{\Delta x}.$$

Since $\rho^{\pm}_{j+1/2} = \alpha(x_{j+1/2})$, the time derivative of momentum is exactly zero. Moreover, the source term and the momentum flux in the energy equation are identically zero at the hydrostatic equilibrium, thus the flux and the source are automatically well balanced, even without the well balanced quadrature of the source (15).

Further, the well balanced reconstruction (12) ensures that $U^{+}_{j+1/2} = U^{-}_{j+1/2} = U_{j+1/2}$ at equilibrium. Thus, for any consistent numerical flux, $\mathcal{F}(U^{+}_{j+1/2}, U^{-}_{j+1/2}) = f(U_{j+1/2})$, which ensures that there is no artificial viscosity. Finally at equilibrium $f(U_{j+1/2}) = [0, p_{j+1/2}, 0] = [0, \beta_{j+1/2}, 0]$. □

**Remark 1.** *The semidiscrete scheme* (17) *is well balanced for* any *consistent numerical flux. We point out that in particular we do not need to assume that the numerical flux is able to preserve contact discontinuities exactly. This allows the application of this framework also to numerical fluxes, such as Lax Friedrichs and the Central Upwind or Rusanov numerical fluxes.*

**First order scheme** Reconstruct the boundary extrapolated data from the cell averages of the auxiliary quantities $\overline{r}, \overline{\pi}, \overline{\rho v}$ with the conservative piecewise constant polynomial. Integrating in time with forward Euler, we obtain a first order accurate well balanced scheme.

**Second order scheme** Reconstruct the boundary extrapolated data from the cell averages of the auxiliary quantities $\overline{r}, \overline{\pi}, \overline{\rho v}$ with a conservative piecewise linear non oscillatory polynomial. Integrating in time with a second order SSP Runge Kutta scheme, we obtain a second order accurate well balanced scheme. For instance, one could use the Minmod limiter to obtain the non oscillatory reconstruction, and Heun's scheme as a Runge Kutta integrator.

**Remark 2.** *Once the semidiscrete scheme is well balanced, thanks to Theorem 1, a Runge-Kutta time integration yields naturally a well balanced fully discrete numerical scheme, because all Runge-Kutta stages will be zero at equilibrium.*

## 3.1 High order accurate well balanced method

High order accuracy can be achieved improving the accuracy of the reconstruction, of the time integration, and of the order of the well balanced quadrature (14). Let $q$ be the desired order of accuracy.

If the reconstruction is computed on the fluctuations from equilibrium, since any high order polynomial reconstruction preserves zero exactly, then at equilibrium, the reconstructed data of the fluctuations will remain zero. Thus, at equilibrium $p^{\pm}_{j+1/2} = \beta(x_{j+1/2})$ and $\rho^{\pm}_{j+1/2} = \alpha(x_{j+1/2})$. So, any high order reconstruction remains well balanced, provided it is applied to the fluctuations (10). To ensure the desired accuracy we will consider a piecewise polynomial reconstruction of order $q$.

Next, one needs to increase the order of accuracy of the time integrator. This, thanks to Theorem 1, does not modify the well balanced property of the scheme. Thus any $q$ order Runge Kutta or multistep scheme can be used to boost the final accuracy of the scheme.

However, the scheme just described is limited to second order accuracy in two key aspects. First of all the cell averages of the pressure, obtained as in (11) are only second order accurate approximations. Secondly, the well balanced quadrature $Q_j$ in (14) is also only second order accurate.

Accurate estimates of the cell averages of the pressure can be obtained in this fashion. First, choose a $q$ order accurate quadrature rule to reconstruct the cell averages of a smooth function $g(x)$:

$$\bar{g}_j = \sum_k w_k g(x_{j,k}) = \frac{1}{\Delta x} \int_{x_{j-1/2}}^{x_{j+1/2}} g(x)\,\mathrm{d}x + O(\Delta x)^q, \tag{18}$$

where $x_{j,k}$ are the quadrature nodes contained in the $j$-th cell $(x_j - 1/2\,h, x_j + 1/2\,h)$.

Next, compute the cell averages of the density fluctuations $\bar{r}$, and reconstruct the point values of the density through (12a) at each quadrature node $x_{j,k}$. Obtain also the point values of momentum $m(x_{j,k})$. Then use the quadrature (18) to compute the cell averages of the kinetic energy

$$\overline{K}_j = \frac{1}{2} \sum_k w_k \frac{m^2(x_{j,k})}{\rho(x_{j,k})}. \tag{19}$$

Note that, since the reconstructions of $\rho$ and $m$ are well balanced, $\overline{K}_j = 0$ at equilibrium. Now, compute the cell averages of the pressure, as $\bar{p}_j = (\gamma - 1)(\overline{E}_j - \overline{K}_j)$. From these, compute the cell averages of the pressure fluctuations $\bar{\pi}_j = \bar{p} - \bar{\beta}$ and reconstruct the point values of the pressure at the cell edges as in (12b). Finally, we obtain the well balanced reconstruction of the energy at the cell edges as

$$E(x_{j+1/2})^\pm = \frac{1}{\gamma - 1} p(x_{j+1/2})^\pm + \frac{1}{2} \frac{(m(x_{j+1/2})^\pm)^2}{\rho(x_{j,k})}. \tag{20}$$

To increase the order of the quadrature rule (16), without loosing the well balanced property, we use the same approach of [18]. The idea is to increase accuracy applying the Richardson extrapolation to the well balanced, second order quadrature rule (14) for momentum, and (15) for the energy equation. This can be done using Romberg's method, see [22, §3.4]. Dropping the index $j$, let $Q_m^{(0)}$ be the well balanced reconstruction $Q$ of (16) applied as a composite rule to the interval $(x_{j-1/2}, x_{j+1/2})$, subdivided into $m = 2^\ell$ equal subintervals of amplitude $\Delta x/m$. Define the recursion

$$Q_m^{(k)} = \frac{2^{2k} Q_m^{(k-1)} - Q_{m/2}^{(k-1)}}{2^{2k} - 1}, \quad k > 0, k \in \mathbb{N}$$

Then the quadrature defined by

$$Q^q = Q_{2^\ell}^{(\ell)} \tag{21}$$

has accuracy $q = 2 + 2\ell$. For example, the formulas

$$Q^4 = \tfrac{1}{3}(4Q_2^{(0)} - Q_1^{(0)}) \qquad q = 4 \qquad (22)$$

$$Q^6 = \tfrac{1}{45}(64Q_4^{(0)} - 20Q_2^{(0)} + Q_1^{(0)}) \qquad q = 6 \qquad (23)$$

are $q = 4$ and $q = 6$ order accurate, respectively.

**Theorem 2.** *The semidiscrete scheme*

$$\frac{\mathrm{d}}{\mathrm{d}t}\overline{U}_j(t) = -\frac{1}{\Delta x}\left[\mathcal{F}\left(U_{j+1/2}^+(t), U_{j+1/2}^-(t)\right) - \mathcal{F}\left(U_{j-1/2}^+(t), U_{j-1/2}^-(t)\right)\right] + Q_j^q(t), \qquad (24)$$

*where $Q_j^q$ is computed in (21), and the boundary extrapolated data are obtained with the well-balanced reconstructions (12) is exactly well balanced on the hydrostatic equilibrium solution $\rho(x,t) = \alpha(x)$ and $p(x,t) = \beta(x)$.*

*Proof.* If the data at time $t$ are in hydrostatic equilibrium, then $\overline{\rho v}$, $\overline{r}_j$ are zero for all $j$, and as a consequence $m(x_{j,k}) = 0$ at each quadrature point, because any piecewise polynomial reconstruction is able to reproduce the zero function exactly. Thus $\overline{K}_j = 0$ in each cell, and the cell averages of the pressure coincide, within the constant $\gamma - 1$, with the equilibrium values of the cell averages of the energy. Thus, $\overline{p}_j = \overline{\beta}_j$ and $\overline{\pi}_j = 0$. Then, the boundary extrapolated data for the pressure are again $p_{j+1/2}^\pm = \beta(x_{j+1/2})$.

The momentum flux reduces to $\mathcal{F}_{j+1/2}^{\rho v} = \beta(x_{j+1/2})$, and the equation for momentum becomes

$$\frac{\mathrm{d}}{\mathrm{d}t}\overline{\rho v}_j(t) = -\frac{\beta(x_{j+1/2}) - \beta(x_{j-1/2})}{\Delta x} + Q_j^q.$$

Since at equilibrium $\rho(x) = \alpha(x)$, for any reconstructed point $x$, the source $Q_j^q$ of (21), becomes a telescopic sum, and reduces to $Q_j^q = -(\beta(x_j + 1/2) - \beta(x_j - 1/2))$, thus balancing momentum exactly.

The remaining part of the proof is identical to the closing argument of the proof of Theorem 1. □

**Remark 3.** *The computation of the cell averages of the kinetic energy and the application of Romberg's algorithm to the estimate of the source requires the evaluation of the reconstruction at several points within each interval. For this reason, it is crucial to adopt a reconstruction algorithm providing* uniform accuracy *within the cell. This can be achieved using the* CWENO *reconstructions of [9] or the improved* CWENOZ *of [8].*

*The* CWENO *algorithm in fact provides the whole reconstruction polynomial on each cell, and not just the value of the reconstructed polynomial at one single point, as in* WENO*. Thus, it is enough to compute the polynomial and then evaluate the result at each point needed. In* WENO *instead, the non linear weights must be recomputed at each reconstruction point, and the weights for points in the interior of the cell may not be positive, or may exist only at the price of reducing accuracy.* ENO *reconstructions [12] would also determine a polynomial with uniform accuracy within a cell, but the stencil required is much larger than in* CWENO*, [9].*

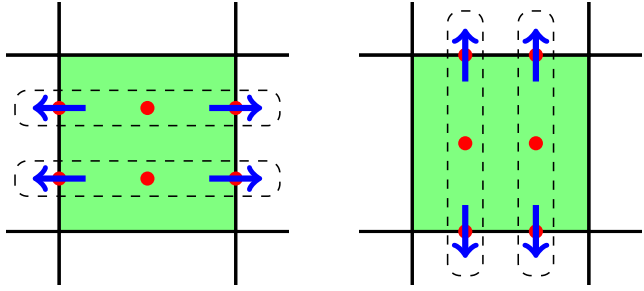

Figure 1: Illustration of the well-balanced quadrature rule for the source term in two space dimensions. Numerical fluxes are represented by the arrows, red dots represent the quadrature nodes for the source terms.

## 4 Extension to higher dimensions

In this section we describe the extension of the proposed numerical method to higher space dimensions. As in the one-dimensional case, we assume that a gravity field $\phi(x)$ is known ($x \in \mathbb{R}^d$) and that a steady-state to be preserved is described by means of a pair of functions $\alpha(x)$ and $\beta(x)$, defined in the entire computational domain, such that $\nabla \beta(x) = \alpha(x) \nabla \Phi(x)$.

We consider a cartesian mesh and compute the cell averages of $\alpha$ and $\beta$ in the preprocessing phase. For simplicity, we describe the algorithm in the case of cells with equal sides, the generalization to rectangular cells being straightforward.

Denoting with $\overline{U}_j$ the cell averages in the $j$-th cell $\Omega_j$, the finite volume form of (1) can be written as

$$\frac{\mathrm{d}}{\mathrm{d}t}\overline{U}_j = -\underbrace{\frac{1}{|\Omega_j|}\int_{\partial\Omega_j} f(U(x,t))\mathrm{d}x}_{K_j} + \underbrace{\frac{1}{|\Omega_j|}\int_{\Omega_j} s(U(x,t))\mathrm{d}x}_{\overline{S}_j}. \qquad (25)$$

To compute $K_j$ a quadrature rule on each face of the cell $\Omega_j$ is needed, while to compute $\overline{S}_j$ one needs a quadrature rule for the volume integral on the cell $\Omega_j$.

In order to minimize the number of numerical flux evaluations, we compute the total flux across each face with a gaussian quadrature rule with weights $\omega_p$ and nodes $\xi_p^k$ belonging to the element $S^k = [-1/2, 1/2]^{d-1}$ orthogonal to the unit vector $e_k$ that represents the $k$-th coordinate direction of $\mathbb{R}^d$. With this notation, in the semidiscrete scheme we have

$$K_j^\Delta = \sum_{k=1}^d \sum_p \omega_p \frac{\mathcal{F}(x_j + \Delta x(1/2 e_k + \xi_p^k)) - \mathcal{F}(x_j - \Delta x(1/2 e_k + \xi_p^k))}{\Delta x}. \qquad (26)$$

The next task is to generalize to the $d$-dimensional case the well-balanced quadrature rule (21). To this end we consider different quadrature rules in each coordinate direction. Consider the pair of opposite faces $x_j \pm \Delta x/2 \, e_k S_k$. The quadrature nodes of (26) on these faces are $x_j + \Delta x(\pm 1/2 e_k + \xi_p)$. Let us denote by $Q^q_{j,k,p}$ the quadrature obtained by applying to the segment with the endpoints $x_j + \Delta x(\pm 1/2 e_k + \xi_p)$ the one-dimensional formula of order $q$ given by Romberg's rule (21). Then the volume average source term in the $k$-th momentum equation is discretized using the cartesian product of the gaussian rule for the face $S_k$ and the one-dimensional well-balanced rule applied in the direction $e_k$, as follows

$$\widetilde{Q}^q_{j,k} = \sum_p \omega_p Q^q_{j,k,p}. \tag{27}$$

For example, for a scheme of order 3 or 4, one employs two gaussian nodes per direction on $S_k$ and the level-2 Romberg formula $Q^4$ of (22). For the case of two space dimensions, the quadrature nodes of $\widetilde{Q}^4_{j,1}$ and, respectively $\widetilde{Q}^4_{j,2}$, are depicted in the left (respectively right) part of Figure 1. The dashed lines in the figure highlight that $\widetilde{Q}^4_{j,k}$'s are cartesian products of the high-order well-balanced quadratures (21) with the gaussian rules on the faces and that each term in the sum (27) is designed to autonomously balance the fluxes at each pair of gaussian quadrature nodes located on opposite faces. The accuracy of $\tilde{Q}_j$ is the smallest between the accuracies of the well balanced quadrature $Q^q_{j,k,p}$, and of the gaussian formula used along the cell boundaries: in the example in the figure, the final accuracy of the complete quadrature is 3.

Finally, the source in the energy equation is the sum of $d$ terms of the form $\rho u_k \partial_{x_k} \Phi$ and each of them is computed with the rule $\widetilde{Q}^q_{j,k}$ applied to (16) with $y = \rho u_k$.

The complete algorithm that computes the right-hand-side of the semi-discrete scheme using the well-balanced reconstructions can be summarized as follows.

1. Compute the cell averages of the variable $r(x)$ as $\bar{r}_j = \bar{\rho}_j - \bar{\alpha}_j$.

2. Compute a reconstruction procedure to obtain point values $r_j(x), (\rho v_k)_j(x)$ for $k = 1, \ldots, d$ from their respective cell averages.

3. Obtain the recontruction of the density as $\rho_j(x) = r_j(x) + \alpha(x)$. Note that this reconstruction coincides with $\alpha(x)$ if the cell averages of the density represent the steady state.

4. Inside each cell, choose a quadrature rule of sufficient order. Compute at the quadrature nodes the point values of density and momentum. From these, compute the cell averages of the kinetic energy $\bar{K}_j$. Use these to find the well balanced cell averages of the pressure fluctuations $\bar{\pi}_j = (\gamma - 1)(\bar{E}_j - \bar{K}_j) - \bar{\beta}_j$.

5. Reconstruct the pressure fluctuations $\pi_j(x)$ from the cell averages $\bar{\pi}_j$ and define the reconstruction of the pressure inside each cell to be $p_j(x) = \pi_j(x) + \beta(x)$. Note that $p_j(x) = \beta(x)$ at equilibrium.

6. Compute the numerical fluxes using the reconstructed values of density, momentum and pressure to find the reconstructed data at the cell boundaries for the conserved quantities.

7. Compute the cell average of the source in the momentum and energy equations using the rules (27).

We point out that the scheme requires the reconstruction of point values at many different locations on the cell boundary and inside the cell and thus it is very efficient to employ a procedure that defines first a reconstruction polynomial everywhere in the cell, that can be later evaluated where needed. For the same reason it is advisable to avoid dimensionally-split reconstructions. In the numerical tests of this work we employ the third order accurate truly two dimensional CWENO reconstruction of [20].

**Theorem 3.** *The semidiscrete scheme*

$$\frac{\mathrm{d}}{\mathrm{d}t}\overline{U}_j(t) = -K_j^\Delta(t) + \overline{S}_j^\Delta(t), \tag{28}$$

*where $K_j^\Delta(t)$ is defined by (26), $\overline{S}_j^\Delta(t)$ is computed with the quadrature $\widetilde{Q}_{j,k}(t)$ defined by (27) along each component and the boundary extrapolated data are computed with the well-balanced reconstructions described in the algorihtm above, is exactly well balanced on the hydrostatic equilibrium solution $\rho(x) = \alpha(x)$ and $p(x) = \beta(x)$.*

*Proof.* The proof descends from the fact that the $\mathbb{R}^d$ scheme is the cartesian product of $d$ one-dimensional well-balanced schemes (see Theorem 2).

Since the reconstruction is well-balanced, the momentum fluxes, at equilibrium, reduce to

$$K_j^{\Delta,\rho v_k} = \frac{1}{\Delta x}\sum_p \omega_p \left(\beta(x_j + \Delta x(1/2 e_k + \xi_p^k)) - \beta(x_j + \Delta x(-1/2 e_k + \xi_p^k))\right).$$

Using again the well-balanced property of the density reconstruction, at each quadrature point $\rho/\alpha = 1$ and thus $\overline{S}_j^{\Delta,\rho v_k}$ is a telescopic sum and reduces to $K_j^{\Delta,\rho v_k}$. □

## 4.1 Well balancing around a moving equilibrium

In this section we describe the variations that must be introduced in the well balanced scheme described above to include also the case in which the equilibrium velocity has a constant non zero component, perpendicular to the gravity field. In particular, we consider the steady solution (9). Since the gravity field has a constant direction, we choose the $y$ axis to be aligned with the gravity field.

The main difficulty here is to obtain well balanced reconstructions of all variables, to ensure that at equilibrium the reconstructed data are continuous across cell interfaces, thus zeroing the artificial diffusion terms in the numerical fluxes, see also §2.

Thus, we start with the well balanced reconstruction of the density, as in (12a). Next, we consider the fluctuations with respect to equilibrium of momentum

$$\overline{\mu}_j = \overline{m}_j - (U\overline{\alpha}_j, 0).$$

Let $\mu_j(x, y)$ be the reconstruction of the momentum fluctuation at the point $(x, y)$ in the cell $\Omega_j$, with the correct order. We reconstruct the point values of momentum as

$$m_j(x, y) = \mu_j(x, y) + (U\alpha(x, y), 0).$$

As before, note that the order of the reconstruction is preserved while subtracting $(U\overline{\alpha}_j, 0)$ and adding $(U\alpha(x, y), 0)$.

The new well balanced reconstruction of momentum is now used to compute the cell average of the kinetic energy, which is needed to obtain the well balanced reconstruction of the pressure. In other words, we find the cell average of the kinetic equation, using the well balanced reconstruction of momentum and density, then we obtain the cell averages of the pressure fluctuations $\overline{\pi}_j = (\gamma - 1)(\overline{E}_j - \overline{K}_j) - \overline{\beta}_j$, from which we compute the well balanced point values of the pressure with (12b).

At equilibrium, $\overline{\mu} \equiv 0$, and thus $m_j(x, y) = (U\alpha(x, y), 0)$, so that the reconstruction is continuous across cell interfaces, guaranteeing a zero artificial diffusion, in all components of the momentum equation. Further, at equilibrium, the direction of $v(x, y)$ is still perpendicular to $\nabla \Phi$, ensuring that the source in the energy equation vanishes. For this reason, it is not necessary to include a correction to the well balanced treatment of the source. Finally, at equilibrium the energy flux $\nabla \cdot (v(E + p))$ becomes $[U, 0] \cdot \nabla^T (E + p)$ which is zero, because the pressure and the total energy are constant in the $x$ direction.

Note that the algorithm described earlier is a special case of the well balanced algorithm presented in this subsection. In fact, if the steady state is at rest, namely $U = 0$, momentum is already an equilibrium variable, or in other terms $\overline{\mu} = \overline{m}$.

## 5  Numerical tests

We start with convergence tests against an exact solution. We choose a variation of the solution proposed in [25], namely

$$\rho(x, t) = 1 + \tfrac{1}{5} \sin(k\pi(x - u_0 t)) \tag{29}$$

$$p(x, t) = \tfrac{9}{2} - (x - u_0 t) + \tfrac{1}{k\pi} \cos(k\pi(x - u_0 t)), \tag{30}$$

where $u_0$ is a constant. Here we take $u_0 = 1$, $\Phi(x) = gx = x$ and $k = 5$, to have a highly oscillatory solution. The Courant number is $c = 0.45$, and the integration is carried out on $[0, 2]$ up to $t_f = 0.1$.

Fig. 2 contains the log-log plots of the error versus the number of grid points $N$. The left panel contains the low order schemes (order 1 in red at the top, and order 2 in blue at the bottom of the figure). The right panel contains the history of convergence for the schemes of order 3 and 5. Note the different vertical scale in the two plots. The different markers represent the errors obtained on density, momentum and energy. Table 1 contains the actual data obtained on the density equation.

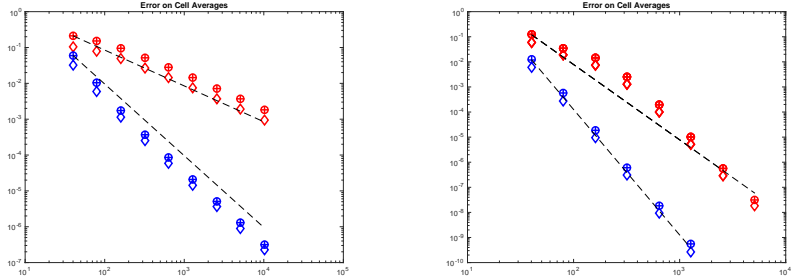

Figure 2: Convergence test for the exact smooth solution of (29) for schemes of order 1 (top), and 2 (bottom curve) of the left figure, and 3 (top), 5 (bottom curve) of the right figure. The errors shown are for density ($\circ$), momentum ($+$) and energy ($\diamond$).

| N | 1 | | 2 | | 3 | | 5 | |
|---|---|---|---|---|---|---|---|---|
| | error | rate | error | rate | error | rate | error | rate |
| 40 | 2.12e-01 | | 5.99e-02 | | 1.22e-01 | | 1.23e-02 | |
| 80 | 1.52e-01 | 0.48 | 1.02e-02 | 2.55 | 3.54e-02 | 1.78 | 5.66e-04 | 4.44 |
| 160 | 9.33e-02 | 0.70 | 1.76e-03 | 2.54 | 1.46e-02 | 1.27 | 1.88e-05 | 4.91 |
| 320 | 5.23e-02 | 0.84 | 3.63e-04 | 2.28 | 2.53e-03 | 2.52 | 5.97e-07 | 4.98 |
| 640 | 2.78e-02 | 0.91 | 8.49e-05 | 2.10 | 1.91e-04 | 3.73 | 1.86e-08 | 5.00 |
| 1280 | 1.43e-02 | 0.95 | 2.08e-05 | 2.03 | 1.04e-05 | 4.20 | 5.69e-10 | 5.03 |
| 2560 | 7.29e-03 | 0.98 | 5.16e-06 | 2.01 | 5.61e-07 | 4.21 | | |
| 5120 | 3.67e-03 | 0.99 | 1.29e-06 | 2.00 | 3.22e-08 | 4.12 | | |
| 10240 | 1.84e-03 | 0.99 | 3.22e-07 | 2.00 | 3.01e-09 | 3.42 | | |

Table 1: Errors on the density and convergence rates against the exact solution (29), for schemes of order 1, 2, 3 and 5 respectively.

## 5.1 Well balanced tests

As in [13], we consider the time needed by sound waves to cross the computational domain

$$\tau_s = 2 \int_a^b \frac{1}{c(x)} \mathrm{d}x, \qquad (31)$$

as a measure to estimate how "large" is the computational time.

|       | $\Phi(x) = x$ | | | $\Phi(x) = x^2$ | | | $\Phi(x) = \sin(2\pi x)$ | | |
|-------|---------|---------|---------|---------|---------|---------|---------|---------|---------|
| order | $\rho$  | $\rho u$ | $E$    | $\rho$  | $\rho u$ | $E$    | $\rho$  | $\rho u$ | $E$    |
| 1 | 6.71e-17 | 1.51e-16 | 3.60e-16 | 8.71e-17 | 1.00e-16 | 3.69e-16 | 3.63e-17 | 1.88e-16 | 3.98e-16 |
| 2 | 1.80e-16 | 1.10e-16 | 3.04e-16 | 3.00e-16 | 1.44e-16 | 3.78e-16 | 1.16e-16 | 1.82e-16 | 3.71e-16 |
| 3 | 3.06e-16 | 1.69e-16 | 5.04e-16 | 3.06e-16 | 1.69e-16 | 5.04e-16 | 5.37e-16 | 3.01e-16 | 7.31e-16 |
| 5 | 4.35e-16 | 2.05e-16 | 6.38e-16 | 5.33e-16 | 2.25e-16 | 7.34e-16 | 1.01e-15 | 3.49e-16 | 1.11e-15 |

Table 2: Distance from the equilibrium solution, for schemes of order 1, 2, 3 and 5 respectively, for isothermal equilibrium.

We start from isothermal equilibrium. We recall that in this case,

$$\alpha(x) = \frac{1}{T} e^{-\frac{\Phi(x)}{T}}, \qquad \beta(x) = e^{-\frac{\Phi(x)}{T}}. \tag{32}$$

The temperature is chosen as $T = 1$. The computational domain is $[0, 1]$. With these choices, $\tau_s \simeq 1.69$. We take $t_f = 2$, so that a typical wave has the time to cross the whole domain. We study the well balancing property for several potentials, $\Phi(x) = x$, $\Phi(x) = x^2$ and $\Phi(x) = \sin(2\pi x)$ on $[0, 1]$. The deviations from equilibrium in the $L^1$ norm appear in Table 2. For the polytropic equilibrium, the results are shown only for the potential $\Phi(x) = x^2$, on the left column of Table 3. For all schemes studied in this work, the equilibrium is maintained within machine precision.

We consider now a well balanced test around an equilibrium solution which is neither isothermal nor polytropic from [5].

$$\alpha(x) = e^{-x}, \qquad \beta(x) = (1+x) e^{-x}. \tag{33}$$

In this case, the temperature is $T(x) = 1 + x$. The system is in equilibrium for the potential $\Phi(x) = -\frac{1}{2} x^2$, and for $u = 0$. So this test is neither polytropic nor isothermal. But since we have an explicit and analytic expression for the density and pressure at equilibrium, it is possible to apply the technique described in this paper, and balance the scheme around the equilibrium (33). We start the computation with initial data coinciding with the equilibrium states given in (33), and integrate as before up to $t_f = 2$, with Dirichlet boundary conditions.

The results for all schemes presented in this work appear in Table 3. Note that all variables are in equilibrium, while the scheme analyzed in [5] is able to guarantee that the velocity remains zero, while in the other variables the error is of the same size of the truncation error, and not machine precision as here.

## 5.2 Perturbations of equilibrium states

The main application of well balanced schemes is in the integration of problems which are small perturbations of equilibrium states. In this case, the failure of a scheme to reproduce steady solutions generates errors large enough to prevent the resolution of the small perturbations. We will show that well-balanced schemes can correctly resolve perturbations as small as the typical size of the local truncation error.

|       | Polytropic, $\Phi(x) = x^2$ |          |          | Non Isothermal |          |          |
|-------|---------|----------|----------|----------|----------|----------|
| order | $\rho$  | $\rho u$ | $E$      | $\rho$   | $\rho u$ | $E$      |
| 1     | 1.63e-16 | 2.29e-16 | 4.43e-16 | 1.17e-17 | 2.75e-16 | 1.94e-16 |
| 2     | 2.38e-16 | 1.61e-16 | 4.66e-16 | 1.95e-16 | 1.93e-16 | 6.57e-16 |
| 3     | 3.80e-16 | 2.20e-16 | 6.78e-16 | 2.30e-16 | 1.68e-16 | 4.76e-16 |
| 5     | 6.06e-16 | 2.39e-16 | 8.78e-16 | 5.22e-16 | 2.41e-16 | 8.70e-16 |

Table 3: Distance from the equilibrium solution, for schemes of order 1, 2, 3 and 5 respectively, for polytropic equilibrium, and for the equilibrium solution of eq. (33).

As an example, we consider a test proposed in [5]. In this test, one starts from an equilibrium isothermal solution of the form (32), with $T = 1$, and potential $\Phi(x) = x^2$, on the interval $[0, 1]$. The data in Table 2 show that all the schemes studied here are exactly well balanced on this solution. As non well balanced schemes, we consider methods in which the reconstruction is carried out on the cell averages of conservative variables $\rho$ and $E$, instead of on fluctuations $r$ and $\pi$ from equilibrium variables.

We run a third order $N3$ and a fifth order $N5$ scheme on the isothermal equilibrium solution (32). These schemes are built with a standard CWENO reconstruction of order 3 and 5 respectively, applied on cell averages of density, momentum and total energy. The source is integrated with the well balanced quadrature. These schemes then are not well balanced, and therefore they give an error with respect to the equilibrium solution, which is of the same size of the local truncation error. For a grid with 40 points on $[0, 1]$, we find a residual with respect to the exact steady state solution which is of order $10^{-4}$ on the velocity and $10^{-4}$ for the pressure in the case of the third order $N3$ scheme, and are $10^{-8}$ on both velocity and pressure for the fifth order scheme $N5$ scheme.

Next, we add a perturbation in the pressure in the initial data. More precisely,

$$\rho(x, t=0) = e^{-\Phi(x)}; \qquad v(x, t=0) = 0; \qquad p(x, t=0) = e^{-\Phi(x)} + Ae^{-100\left(x - \frac{1}{2}\right)^2},$$

where we can modify the amplitude of the perturbation $A$ with respect to the equilibrium solution. For $A = 10^{-3}$, the perturbation has a size close to the truncation error we noticed above for the third order scheme. Running all schemes on this problem, with $N = 120$ and $t_f = 0.25$, we find the pressure and velocity fluctuations of Fig. 3. The figure shows that all schemes resolve the perturbations, albeit the third order, not well-balanced scheme $N3$ is losing accuracy with respect to the well balanced schemes. On this scale, the truncation error of the $N5$ scheme is so small that it is able to resolve the perturbation even though it is not well balanced.

Decreasing the amplitude of the perturbation to $A = 10^{-5}$, the third order

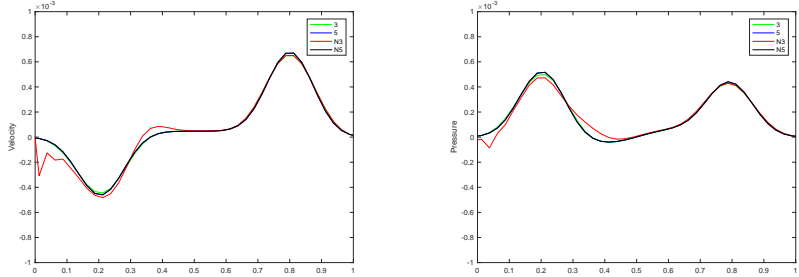

Figure 3: Perturbation of an equilibrium solution of amplitude $A = 10^{-3}$ for well balanced schemes of order 3 and 5 (green and blue curves), and non well balanced schemes of order 3 and 5 (red and black curves). $N = 40$ points per unit length.

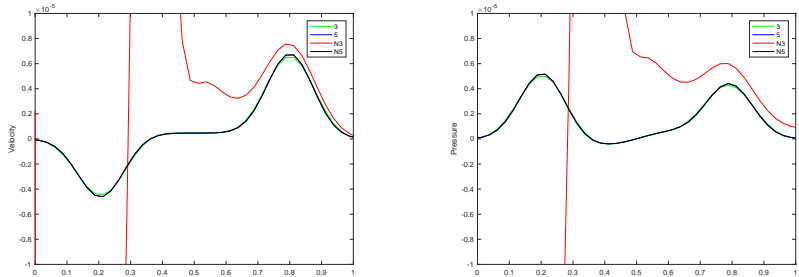

Figure 4: Perturbation of an equilibrium solution of amplitude $A = 10^{-5}$ for well balanced schemes of order 3 and 5 (green and blue curves), and non well balanced schemes of order 3 and 5 (red and black curves). $N = 40$ points per unit length.

$N3$ scheme is no longer able to resolve the perturbation, which is too small with respect to its residual on the equilibrium solution, red curve in Fig 4. Both well balanced schemes, of order 3 and 5, instead resolve the small perturbation, as does the fifth order $N5$.

Finally, decreasing further the amplitude of the perturbation to $A = 10^{-7}$, we obtain the results of Fig. 5. The third order $N3$ misses the solution completely, and it is not shown. In this case, the fifth order $N5$ scheme is also starting to lose resolution, because now the magnitude of the perturbation is reaching the magnitude of its error on the equilibrium solution. For this test, we also show a reference solution (dashed line) obtained with the fifth order well balanced scheme. It is noteworthy that the other schemes obtain a comparable precision with only 40 points. The results on the perturbation of an equilibrium solution are shown on the interval $[0, 1]$, but they were computed on the larger interval $[-1, 2]$, to avoid perturbations coming from inexact boundary conditions, which do not have the time to reach the region shown.

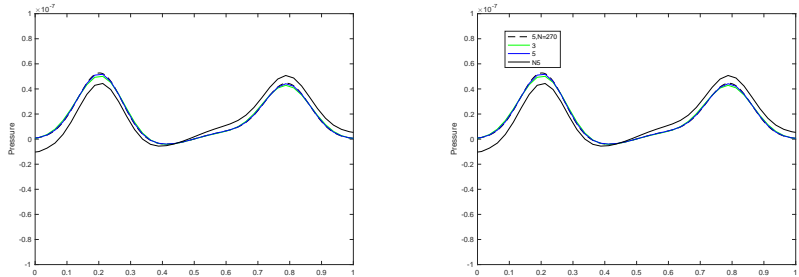

Figure 5: Perturbation of an equilibrium solution of amplitude $A = 10^{-7}$ for well balanced schemes of order 3 and 5 (green and blue curves), and non well balanced scheme of order 5 (black curve). $N = 40$ points per unit length. The black, dashed curve is obtained with the well balanced scheme of order 5 and 270 points.

## 5.3 Accuracy of the 2D scheme

All results shown in this section were obtained with the third order well balanced scheme.

To test the accuracy of the third order scheme, we consider the two-dimensional generalization of (29), namely the exact solution defined by

$$\begin{aligned} \rho(x,y,t) &= 1 + \tfrac{1}{5} \sin(k\pi((x+y) - (u_0+v_0)t)) \\ p(x,y,t) &= \tfrac{9}{2} - (x+y - (u_0+v_0)t) + \tfrac{1}{5k\pi} \cos(k\pi(x+y - (u_0+v_0)t)), \end{aligned} \quad (34)$$

with gravity potential $\Phi(x, y) = x + y$ and parameter $k = 1$. The computational domain is $[0, 2]^2$ and the final time is $t = 0.1$.

This is not an equilibrium solution, and we use it to test the accuracy of the well-balanced scheme; in particular we choose $\alpha$ and $\beta$ corresponding to an isothermal problem (6), with temperature $T_{eq}$ obtained from the average of the temperature of the initial condition.

In Table 4 we report the 1-norm of the errors obtained for this test on several $N \times N$ grids and the experimental convergence rates. The scheme matches the expected third order accuracy for all components. The table shows results only for density and energy.

## 5.4 Well-balancing of the 2D scheme

We now present a series of tests demonstrating the well-balancing properties of the two-dimensional numerical scheme.

|      | density          |       | energy           |       |
|------|------------------|-------|------------------|-------|
| N    | error            | rate  | error            | rate  |
| 20   | 6.83e-03         |       | 8.91e-03         |       |
| 40   | 8.61e-04         | 2.99  | 1.16e-03         | 2.94  |
| 80   | 1.08e-04         | 3.00  | 1.46e-04         | 2.99  |
| 160  | 1.35e-05         | 3.00  | 1.82e-05         | 3.00  |
| 320  | 1.68e-06         | 3.00  | 2.27e-06         | 3.00  |
| 640  | 2.10e-07         | 3.00  | 2.84e-07         | 3.00  |
| 1280 | 2.63e-08         | 3.00  | 3.55e-08         | 3.00  |

Table 4: 1-norm of error on $N \times N$ grid and convergence rates for test (34)

In the first test we consider the gravity field $\Phi(x,y) = x+y$ and an initial data which is an isothermal steady-state with $T_{eq} = 1/1.21$ as in [25]:

$$\rho_0^{\text{iso}}(x,y) = \frac{e^{-\Phi(x,y)/T_{eq}}}{T_{eq}} \qquad p_0^{\text{iso}}(x,y) = e^{-\Phi(x,y)/T_{eq}}. \qquad (35)$$

In this test, $\alpha$ and $\beta$ are chosen to enforce the isothermal equilibrium from (6), therefore $\alpha = \rho^{\text{iso}}$ and $\beta = p^{\text{iso}}$.

Next, we consider an hydrostatic polytropic steady state defined by

$$\rho_0^{\text{poly}}(x,y) = \left(1 - \tfrac{\nu-1}{\nu}\Phi(x,y)\right)^{\frac{1}{\nu-1}} \qquad p_0^{\text{poly}}(x,y) = (\rho(x,y))^{\nu} \qquad (36)$$

and in particular the case with $\nu = 1.2$. Here, we choose $\alpha$ and $\beta$ corresponding to polytropic equilibrium, from equation (7).

Finally, we show results obtained on an equilibrium which is neither isothermal nor polytropic. In particular we consider the gravity potential $\Phi(x,y) = r^2$, where $r = \sqrt{(x^2+y^2)}$ and the initial data

$$\rho^{\text{gen}}(r,\theta) = e^{-r} \qquad p^{\text{gen}}(r,\theta) = (1+r)e^{-r}. \qquad (37)$$

This is a steady-state with equilibrium temperature $T(x,y) = 1 + r$. It is a radial version of a test presented in one space dimension in [5]. Accordingly, the functions $\alpha$ and $\beta$ are chosen as in (8).

Therefore in all cases, $\alpha(x,y) = \rho(x,y,t=0)$ and $\beta(x,y) = p(x,y,t=0)$. The initial data are evolved up to $t_f = 0.1$. We report the deviations from the equilibrium solution in Table 5. We see that in all three tests the data are of the order of the machine precision. Thus the scheme is exactly well-balanced.

## 5.5 Perturbation of an isothermal steady state

In this test we consider the gravity field $\Phi(x,y) = x+y$ and an initial data which is a perturbation of the isothermal steady-state (35):

$$\rho_0 = \rho_0^{\text{iso}} \quad u = v = 0 \quad p_0 = p_0^{\text{iso}} + Ae^{-100((x-0.3)^2+(y-0.3)^2)/T_{eq}}.$$

Isothermal

| N | density | momX | momY | energy |
|---|---|---|---|---|
| 20 | 2.80e-17 | 3.51e-17 | 3.70e-17 | 1.60e-16 |
| 40 | 3.14e-17 | 7.63e-17 | 6.98e-17 | 2.16e-16 |
| 80 | 2.78e-17 | 8.36e-17 | 8.62e-17 | 2.10e-16 |
| 160 | 3.90e-17 | 2.16e-16 | 2.12e-16 | 3.38e-16 |
| 320 | 4.65e-17 | 3.36e-16 | 3.20e-16 | 3.93e-16 |
| 640 | 8.76e-17 | 8.70e-16 | 8.57e-16 | 1.20e-15 |
| 1280 | 9.30e-17 | 1.13e-15 | 1.14e-15 | 1.21e-15 |

Polytropic

| N | density | momX | momY | energy |
|---|---|---|---|---|
| 20 | 5.50e-17 | 6.77e-17 | 6.62e-17 | 1.92e-16 |
| 40 | 8.35e-17 | 1.05e-16 | 1.03e-16 | 2.78e-16 |
| 80 | 1.07e-16 | 1.13e-16 | 1.10e-16 | 3.33e-16 |
| 160 | 9.97e-17 | 2.90e-16 | 2.91e-16 | 3.31e-16 |
| 320 | 6.62e-17 | 4.99e-16 | 4.99e-16 | 3.06e-16 |
| 640 | 9.88e-17 | 1.17e-15 | 1.16e-15 | 6.25e-16 |
| 1280 | 8.54e-17 | 1.72e-15 | 1.72e-15 | 7.91e-16 |

Neither isothermal nor polytropic

| N | density | momX | momY | energy |
|---|---|---|---|---|
| 20 | 0.00e+00 | 1.70e-16 | 1.74e-16 | 5.66e-16 |
| 40 | 1.11e-18 | 2.08e-16 | 2.05e-16 | 6.77e-16 |
| 80 | 4.15e-18 | 3.50e-16 | 3.41e-16 | 8.42e-16 |
| 160 | 7.49e-18 | 4.10e-16 | 3.67e-16 | 8.97e-16 |
| 320 | 1.88e-17 | 1.03e-15 | 9.97e-16 | 1.14e-15 |
| 640 | 2.06e-17 | 1.29e-15 | 1.16e-15 | 1.12e-15 |
| 1280 | 5.19e-17 | 3.63e-15 | 3.48e-15 | 2.48e-15 |

Table 5: Errors in 1-Norm on equilibrium solutions: isothermal (35) on the top, polytropic (36) in the middle and the equilibrium (37) on the bottom.

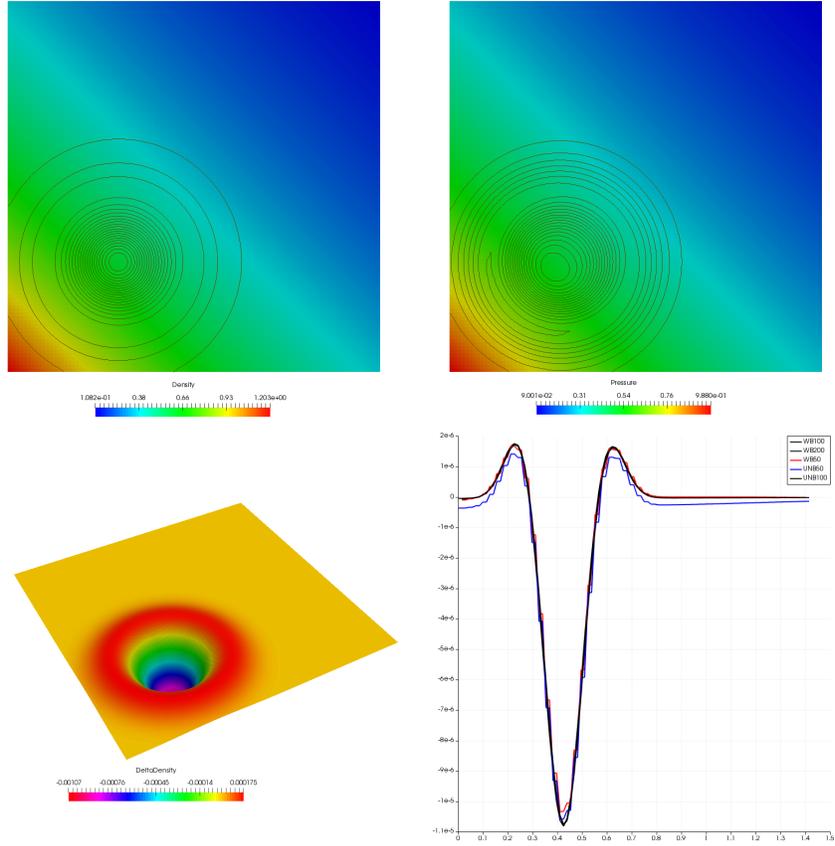

Figure 6: Top: density (left) and pressure (right). The contours represent the perturbations (20 isolines between -0.001 and 0.0002) for $\rho$ and 20 isolines between -0.00027 and 0.00027 for $p$). Bottom left: extrusion of the density perturbation. Bottom right: comparison of the density for different grids and schemes (black lines for $\Delta x = 1/100$ and $\Delta x = 1/200$ balanced, $\Delta x = 1/100$ WB; red for $\Delta x = 1/50$ WB; blue for $\Delta x = 1/50$ UNB)

In particular, we discuss the results for $A = 10^{-5}$. From Fig. 4 we have seen that for a perturbation of a steady state solution with amplitude $10^{-5}$ the unbalanced 1D scheme of order 3 was unable to detect the behaviour of the solution, while the well balanced version is able to compute an accurate solution.

The numerical simulations with the well-balanced (WB) scheme are performed in the domain $[0, 1] \times [0, 1]$. They are compared with computations by a not well-balanced (UNB) scheme in which the reconstructions are applied directly to the conservative variables. In the scheme UNB, disturbances due to the boundary conditions are clearly visible, but in order to make a fairer comparison of the schemes, we have employed the larger domain $[-1, 2] \times [-1, 2]$ for the unbalanced simulations.

The density and pressure fields at final time $t = 0.15$ obtained with well-balanced schemes are presented in the top row of Figure 6. In the same plots, the contours of the density and pressure perturbations are shown. The density perturbation is also depicted as an extruded surface in the bottom-left panel of the same figure. The solutions with different schemes are compared in the bottom-right panel of the figure.

We use the solution obtained by the WB scheme with $\Delta x = 1/200$ as reference solution. The solution on the $\Delta x = 1/100$ grid with WB is almost coincident with the reference solution, showing convergence. At this resolution also the unbalanced UNB scheme captures the correct behaviour. The solutions computed on $\Delta x = 1/50$ grids are different from each other. The one computed with the WB scheme is very close to the reference, even on such a coarse grid. On the other hand, the solution computed with the UNB scheme is quite far away from the others. The unbalanced scheme, in fact, is not able to preserve the isothermal background (see the differences at the extreme left and right of the section) and as a consequence also computes a less accurate solution inside the perturbation.

## 5.6 Radial Rayleigh-Taylor instability

Let $(r, \theta)$ denote the polar coordinates and consider the radial gravity potential $\Phi(r, \theta) = r$. In this case, the solution $\rho = p = e^{-r}$ is an isothermal equilibrium. We consider the following perturbed initial data:

$$\rho(r, \theta) = \begin{cases} e^{-r} & r < r_0 \\ \frac{1}{a} e^{\frac{-r + r_0(1-a)}{a}} & r > r_0 \end{cases} \qquad p(r, \theta) = \begin{cases} e^{-r} & r < r_I(\theta) \\ e^{\frac{-r + r_0(1-a)}{a}} & r > r_I(\theta). \end{cases}$$

Choosing $a = \frac{e^{-r_0}}{e^{-r_0} + \Delta \rho}$, the pressure is continuous, but a jump discontinuity of size $\Delta \rho$ is present on the interface defined by $r = r_I(\theta) = r_0(1.0 + \eta \cos(k\theta))$. These initial data give rise to a Rayleigh-Taylor instability, whose onset is favoured by the wiggles in the discontinuity line $r_I(\theta)$.

In Figure 7 we show the time evolution computed with our third order well-balanced scheme for $r_0 = 0.5$, $k = 20$, $\eta = 0.02$ and $\Delta \rho = 0.1$, following [15].

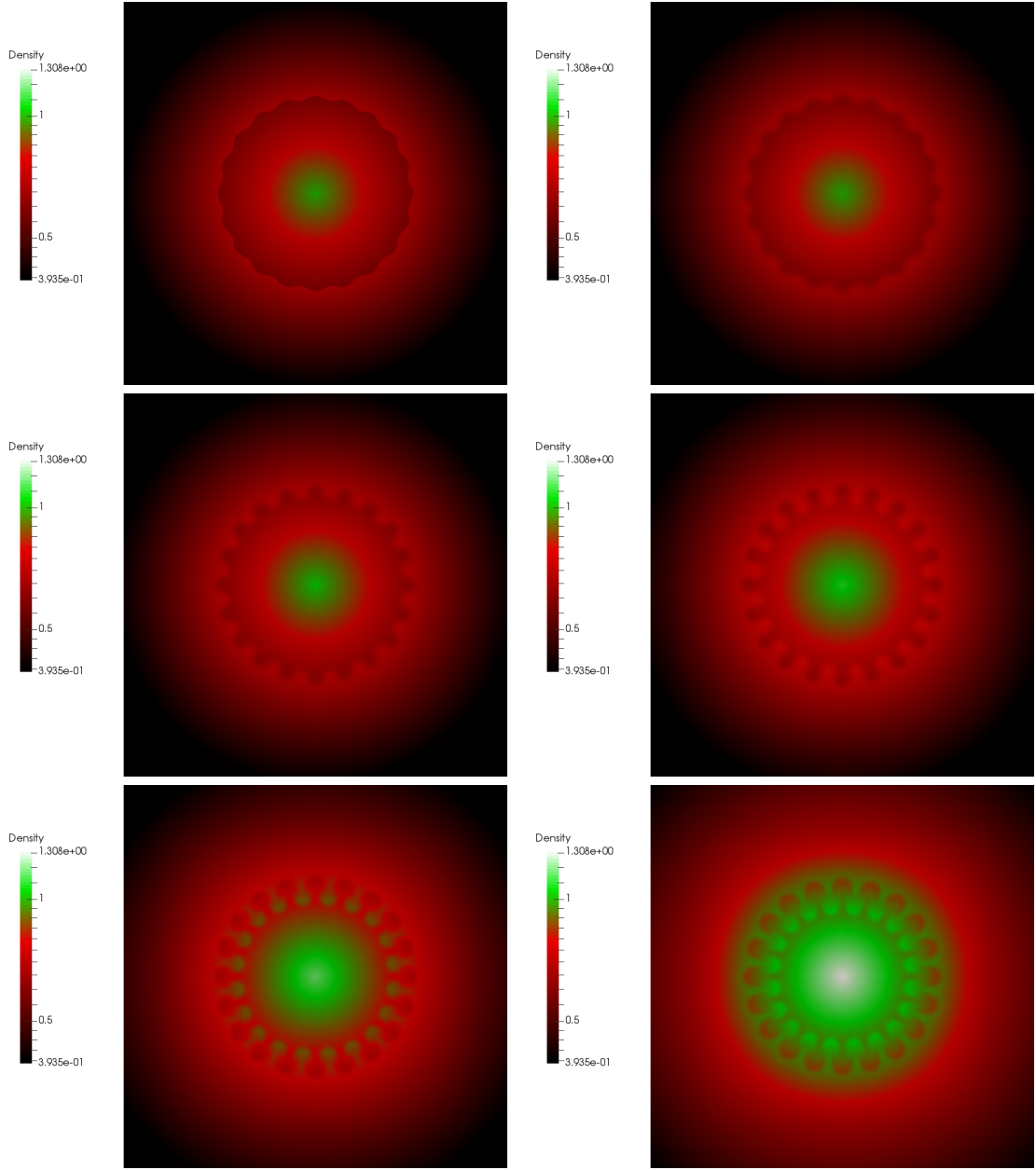

Figure 7: Third order well balanced scheme on $400 \times 400$ grid: times 0 to 4.0 equispaced

|     | density      |      | energy       |      |
| --- | ------------ | ---- | ------------ | ---- |
| N   | error        | rate | error        | rate |
| 20  | 6.83e-03     | -    | 8.93e-03     | -    |
| 40  | 8.62e-04     | 2.98 | 1.16e-03     | 2.94 |
| 80  | 1.08e-04     | 2.99 | 1.46e-04     | 3.00 |
| 160 | 1.35e-05     | 3.01 | 1.82e-05     | 3.00 |
| 320 | 1.69e-06     | 2.99 | 2.27e-06     | 3.00 |

Table 6: Accuracy on moving equilibrium.

The perturbation of the steady state is localized around the initial position of the interface. Plume like structures form from the initial wiggles in $r_I(\theta)$. This level of accuracy is only possible with a scheme that is well-balanced and thus can preserve exactly the equilibrium solution away from the interface, up to the boundary of the computational domain. Note also the lack of grid orientation artifacts: such spurious effects do not occur in our scheme, thanks to the high order accuracy of the scheme *and* to the truly 2D reconstruction of [20]. We further point out that unbalanced schemes are not only less accurate (see the previous test), but also typically generate spurious signals at the boundary of the computational domain that would overwhelm the weak perturbation of the isothermal background equilibrium.

## 5.7 Numerical tests on moving equilibrium

In this section we consider the numerical scheme described in §4.1, designed to preserve moving gas equilibria, with a constant speed perpendicular to $\nabla\Phi$. First we verify that the new reconstruction does not pollute the order of accuracy. We consider the data leading to the exact solution (34). The discretization parameters are the same as in §5.3, while the constant horizontal velocity is $U = 1$. We exhibit the results in Table 6 for density and energy, where the third order accuracy is apparent. Note that in this test the gravity field is skew with respect to the grid and to the imposed constant velocity field.

Next, we test the well balancing property of the scheme around the moving isothermal equilibrium with density and pressure prescribed in (35), gravity field $\Phi(x, y) = y$, and constant velocity field $(u_0, v_0) = (1, 0)$. The errors found in this case are reported in the top part of Table 7, from which it is clear that the novel scheme is indeed well balanced. The lower part of the table reports the errors obtained with a scheme that is well-balanced with respect to the stationary isothermal equilibrium. It is clear that the errors on the horizontal momentum decrease with the grid size, but the scheme is not well-balanced in this case.

| $N$ | density | momX | momY | energy |
|-----|---------|------|------|--------|
| 20  | 5.39e-13 | 5.39e-13 | 1.09e-12 | 2.88e-12 |
| 40  | 3.16e-14 | 3.16e-14 | 6.32e-14 | 1.95e-13 |
| 80  | 1.87e-15 | 1.88e-15 | 3.89e-15 | 1.26e-14 |
| 160 | 5.59e-17 | 6.07e-17 | 3.20e-16 | 5.98e-16 |
| $N$ | density | momX | momY | energy |
| 20  | 3.72e-10 | 3.42e-05 | 1.59e-10 | 3.42e-05 |
| 40  | 3.00e-11 | 7.66e-06 | 9.02e-12 | 7.66e-06 |
| 80  | 2.37e-12 | 1.68e-06 | 5.30e-13 | 1.69e-06 |
| 160 | 1.68e-13 | 3.68e-07 | 3.57e-14 | 3.68e-07 |

Table 7: Well balancing errors for moving equilibrium. Top: well-balanced scheme with respect to the moving isothermal equilibrium. Bottom: well-balanced scheme with respect to the stationary isothermal equilibrium.

## 6 Conclusions

In this paper, we have presented high order numerical schemes, which are exactly well balanced around a prescribed equilibrium solution, known through two functions $\alpha$ and $\beta$ such that $\nabla\beta = \alpha\nabla\Phi$. The well balanced property has been proven for schemes of arbitrary order of accuracy and space dimensions. The results have been tested for schemes up to fifth order in 1D, and for third order in 2D, but the method is clearly explained for any order and for any number of space dimensions.

Further, for the particular case in which the gravity field has a constant direction, the scheme can also balance exactly equilibrium flows, with a non zero constant component of the velocity perpendicular to the direction of the the gravitational field. We believe that this is the first high order, i.e. higher than 2, well balanced scheme for this kind of multidimensional equilibrium solution.

The scheme is based on two simple ingredients. First, the reconstruction of point values for the evaluation of the numerical fluxes is performed on *fluctuations* from the equilibrium solution for the density and the pressure, plus the moementum, in the case of moving equilibria. Next, the reconstruction of conservative variables is carried out, respecting the equilibrium information contained in the well balanced pressure and density. Secondly, as in [18], we observe that a well balanced cell average of the source can be easily obtained at second order, and then accuracy can be boosted with extrapolation techniques as Romberg's.

The computation of the quadrature of the numerical fluxes and of the source requires to estimate the numerical solution at several quadrature points within each cell, as is always the case for high order finite volume schemes for balance laws, even without the complications due to the well balanced quadrature for the source. The computational complexity for a very high order scheme can be controlled using reconstructions that provide a single polynomial which is uniformly accurate throughout the cell, as in the case of the CWENO reconstructions of [9], which have been used for the one-dimensional scheme, while the 2D tests were performed with the truly multidimensional third order accurate CWENO reconstruction of [20].

At present, there are very few well balanced schemes for Euler equations with gravity, that are not restricted to a particular order of accuracy, and none for moving equilibria. The class of schemes proposed here is characterized by enforcing equilibrium, without any restriction on the numerical flux function, which needs only to be consistent. This enables to use also very popular numerical fluxes, such as Lax Friedrichs' and its variants, or numerical fluxes which are tailored to particular needs.

We also observe that the schemes proposed here are based on the knowledge of the analytic expression of the equilibrium functions $\alpha$ and $\beta$, but these are not really necessary. As the scheme is presently built, it has a fixed point which coincides exactly with the cell averages of $\alpha$ and $\beta$. However, it is also possible to initialize the scheme without the analytical expression of the equilibrium functions, but with a sufficiently accurate knowledge of their cell averages. We are currently exploring automatic well balancing along these lines.

### Acknowledgments

We would like to thank Praveen Chandrashekar for his precious suggestions. This work was partially funded by DAAD-MIUR Joint mobility program 2015, and INDAM-GNCS 2017 project.# References